\documentclass[amsmath,amssymb,floatfix]{amsart}
\usepackage{graphicx}
\usepackage{dcolumn}
\usepackage{bm}
\usepackage{amssymb}
\usepackage{epsfig}
\usepackage{color}
\usepackage{amsthm}
\usepackage{amsmath}

\newcommand{\ba}{\begin{eqnarray}}
\newcommand{\ea}{\end{eqnarray}}
\newcommand{\be}{\begin{equation}}
\newcommand{\ee}{\end{equation}}
\newcommand{\bdisplay}{\begin{displaymath}}
\newcommand{\edisplay}{\end{displaymath}}

\newtheorem{thm}{Theorem}

\newcommand{\eq}[1]{Eq.\,(\ref{#1})}

\makeatletter
\def\eqnarray{\stepcounter{equation}\let\@currentlabel=\theequation
\global\@eqnswtrue
\tabskip\@centering\let\\=\@eqncr
$$\halign to \displaywidth\bgroup\hfil\global\@eqcnt\z@
  $\displaystyle\tabskip\z@{##}$&\global\@eqcnt\@ne
  \hfil$\displaystyle{{}##{}}$\hfil
  &\global\@eqcnt\tw@ $\displaystyle{##}$\hfil
  \tabskip\@centering&\llap{##}\tabskip\z@\cr}

\def\endeqnarray{\@@eqncr\egroup
      \global\advance\c@equation\m@ne$$\global\@ignoretrue}

\def\@yeqncr{\@ifnextchar [{\@xeqncr}{\@xeqncr[5pt]}}
\makeatother

\begin{document}

\title{Complex asymptotics in $\lambda$ for the Gegenbauer \\ functions $C_\lambda^\alpha(z)$ and $D_\lambda^\alpha(z)$  with $z\in(-1,1)$ }

\author{Loyal Durand }

\address{ Department of Physics, University of Wisconsin-Madison, Madison, WI 53706}

\email{ldurandiii@comcast.net}

\thanks{Current address: 415 Pearl Court, Aspen, CO 81611}

%\note{Current address: 415 Pearl Court, Aspen, CO 81611; ldurandiii@comcast.net} 

\begin{abstract}
We derive asymptotic results for the Gegenbauer functions  $C_\lambda^\alpha(z)$ and $D_\lambda^\alpha(z)$ of the first and second kind for complex $z$ and the degree $|\lambda \lvert\rightarrow\infty$, apply the results to the case  $z\in (-1,1)$, and  establish the connection of these results to asymptotic Bessel-function approximations of the functions for $z\rightarrow \pm 1$.
\end{abstract}

\maketitle

%\keyword{Gegenbauer functions; asymptotics }

%\MSC{33C45,  33C20, 34L10, 30E20}

%%%%%%%%%%%%%%%%%%%%%%%%%
%%%%%%%%%%%%%%%%%%%%%%%%%

\section{Introduction \label{sec:intro}}

A number of results are known for the asymptotic behavior of the Gegenbauer functions of the first and second kinds, $C_\lambda^\alpha(z)$ and $D_\lambda^\alpha(z)$, for the degree $\lvert\lambda \rvert\rightarrow\infty$ with $z$ complex \cite[Sec. 2.3.2\,(17)]{HTF}, \cite[Sec. 6]{DFS}, \cite[Appendix]{LDaddition}, \cite[Sec. 2.3\,(1)]{cohl}, \cite[Sec. 14]{DLMF}. However, as usually stated, the simple results presented here in Thms.~1 and 2 exclude the important cases with $z$ real, $-1<z<1$ and $1<z<\infty$, and exclude the limits $z\rightarrow \pm1$. These cases have been of interest in recent problems, for example, in \cite{cohl}.\footnote{Private communication from Dr.\ Howard Cohl.} 

In the present work, we sketch the derivation of these results, show that they can be extended to include the cases usually excluded. We show also that  the results connect smoothly for $\lvert\lambda \rvert$ large and $\lvert 1\pm z\rvert$ small to asymptotic expansions for $C_\lambda^\alpha(z)$ and $D_\lambda^\alpha(z)$ in terms of Bessel functions, Thms.~3 and 4. Those expansions include the limits $z\rightarrow\pm 1$.

%%%%%%%%%%%%%%%%%%%%%
%%%%%%%%%%%%%%%%%%%%%

\section{Asymptotic results for $C_\lambda^\alpha(z)$ and $D_\lambda^\alpha(z)$ \label{sec:asymptotics}}

%
% Thm. 1 %
\begin{thm}
Let $z\in {\mathbb C}$ and define $z_\pm=z\pm\sqrt{z^2-1}$ with $-\pi\leq{\rm arg}(z\pm 1)\leq\pi$.  Then for  $\Re\lambda\geq 0$ , $\Re\alpha>0$, $-\pi/2\leq {\rm arg}\lambda\leq\pi/2$, and $\lvert\lambda\rvert\rightarrow\infty$,
\ba
\label{thm1D}
D_\lambda^\alpha(z) &=& e^{i\pi\alpha} \frac{ 2^{-\alpha}}{\Gamma(\alpha)}\lambda^{\alpha-1}(z^2-1)^{-\alpha/2}z_+^{-\lambda-\alpha} \left[1+{\mathcal O}(1/\lvert\lambda\rvert)\right],  \\
\label{thm1C}
C_\lambda^\alpha(z) &=& \frac{2^{-\alpha}}{\Gamma(\alpha)}\lambda^{\alpha-1}(z^2-1)^{-\alpha/2}\left(e^{i\pi\alpha} z_+^{-\lambda-\alpha} 
+  z_-^{-\lambda-\alpha} \right)\left[1+{\mathcal O}(1/\lvert\lambda\rvert)\right].
 \ea
\end{thm}

\begin{thm}
% Thm. 2 %
For $z= x$ real with $x\in(-1,1)$, define $x=\cos{\theta}$, $0<\theta<\pi$. Then for  $\lvert\lambda\rvert\rightarrow\infty$ with   $\Re\lambda\geq 0$, $\Re\alpha>0$, $-\pi/2\leq {\rm arg}\lambda\leq\pi/2$, and  $\sqrt{1-x^2}=\sin{\theta}\gg 1/\lvert\lambda\rvert$, the Gegenbauer functions ${\mathsf D}_\lambda^\alpha(\cos{\theta})$ and ${\mathsf C}_\lambda^\alpha(\cos{\theta})=C_\lambda^\alpha(\cos{\theta})$ ``on the cut''  \cite[7, 8]{Askey} have the limiting  behavior
\ba
\label{thm2D}
{\mathsf D}_\lambda^\alpha(\cos{\theta}) &= &  
-\frac{2^{-\alpha+1}}{\Gamma(\alpha)}\lambda^{\alpha-1}(\sin{\theta})^{-\alpha}\sin{((\lambda+\alpha)\theta-\pi\alpha/2)}\left[1+{\mathcal O}(1/\lvert\lambda\rvert)\right], \\
\label{thm2C}
{\mathsf C}_\lambda^\alpha(\cos{\theta}) &=& 
 \frac{2^{-\alpha+1}}{\Gamma(\alpha)}\lambda^{\alpha-1}(\sin{\theta})^{-\alpha}\cos{\left((\lambda+\alpha)\theta-\pi\alpha/2\right)}\left[1+{\mathcal O}(1/\lvert\lambda\rvert)\right].
\ea
These functions are proportional to the Ferrers functions of \cite[Sec. 14.23]{DLMF}).
\end{thm}

\begin{thm}
% Thm. 3 %
For $z$  complex with $z\approx 1$, define $Z=\sqrt{2(\lambda+\alpha)^2(1-z)}$ and $Z'=\sqrt{2(\lambda+\alpha)^2(z-1)}$. Then for $\Re(\lambda+\alpha)\geq 0$, $\Re\alpha\geq-\frac{1}{2}$,    $\lvert \sqrt{z-1}\rvert\ll 1/\lvert\lambda\rvert^{1/3}$,  and $\lvert\lambda\rvert\rightarrow\infty$,
\ba
\label{thm3D}
D_\lambda^\alpha(z) &=&  \frac{1}{\sqrt{\pi}}e^{i\pi\alpha}\frac{1}{\Gamma(\alpha)}2^{-\alpha+\frac{1}{2}}(\lambda+\alpha)^{\alpha-\frac{1}{2}}(z^2-1)^{-\frac{\alpha}{2}+\frac{1}{4}}K_{\alpha-\frac{1}{2}}(Z') \\
&& \times \left[1+{\mathcal O}\left(1/\lvert\lambda\rvert^{2/3}\right)\right], \nonumber \\
C_\lambda^\alpha(z) &=& \frac{\sqrt{\pi}}{\Gamma(\alpha)}2^{-2\alpha+1}(\lambda+\alpha)^{2\alpha-1}\left(\frac{Z}{2}\right)^{-\alpha+\frac{1}{2}}J_{\alpha-\frac{1}{2}}(Z) \\
&& \times\left[1+{\mathcal O}\left(1/\lvert\lambda\rvert^{2/3}\right)\right]. \nonumber
\ea
\\
For $x=\cos{\theta}\in(-1,1)$, $\theta\approx 0$, $X=\sqrt{2(\lambda+\alpha)^2(1-x)}$, and $\lvert\lambda\rvert\rightarrow\infty$, 
\ba
\label{Dcut2}
  {\mathsf D}_\lambda^\alpha(x) &=& -\frac{\sqrt{\pi}}{\Gamma(\alpha)}2^{-\alpha+\frac{1}{2}}(\lambda+\alpha)^{\alpha-\frac{1}{2}}(1-x^2)^{-\alpha/2+\frac{1}{4}}Y_{\alpha-\frac{1}{2}}(X) \\
  && \times\left[1+{\mathcal O}\left(1/\lvert\lambda\rvert^{2/3}\right)\right], \\ \nonumber
\label{Ccut2}
{\mathsf C}_\lambda^\alpha(x) &=& \frac{\sqrt{\pi} }{\Gamma(\alpha)}2^{-\alpha+\frac{1}{2}}(\lambda+\alpha)^{\alpha-\frac{1}{2}}(1-x^2)^{-\alpha/2+\frac{1}{4}}J_{\alpha-\frac{1}{2}}(X) \\
&& \times \left[1+{\mathcal O}\left(1/\lvert\lambda\rvert^{2/3}\right)\right] . \nonumber
\ea
The results in Thms.\ 2 and 3 match in their common range of validity, $1/\lvert\lambda\rvert\ll\lvert \sqrt{z-1}\rvert\ll 1/\lvert\lambda\rvert^{1/3}$.
\end{thm}

\begin{thm}
% Thm. 4 %
For $z$  complex with $z\approx -1$, define $Z''=\sqrt{2(\lambda+\alpha)^2(1+z)}$. Then for $\Re(\lambda+\alpha)\geq 0$, $\Re\alpha\geq-\frac{1}{2}$,    $\lvert \sqrt{z+1}\rvert\ll 1/\lvert\lambda\rvert^{1/3}$,  and $\lvert\lambda\rvert\rightarrow\infty$,
\ba
D_\lambda^\alpha(z) &=& e^{i\pi\alpha}2^{-\alpha}(\lambda+\alpha)^{\alpha-1}\left(2(1+z)\right)^{-\alpha}\frac{1}{\Gamma(\alpha)}\frac{\sqrt{\pi}}{\sin{\pi(\alpha-\frac{1}{2})}}\left(\frac{Z''}{2}\right)^{\frac{1}{2}}  \\
\label{thm4D}
&& \times e^{\mp i\pi(\lambda+2\alpha)}\left[-J_{\alpha-\frac{1}{2}}(Z'')+e^{\pm i\pi(\alpha-\frac{1}{2}}J_{-\alpha+\frac{1}{2}}(Z'')\right] \nonumber \\ 
&& \times\left[1+{\mathcal O}\left(1/\lvert\lambda\rvert^{2/3}\right)\right], \nonumber
\ea
where the + and - signs hold for $z$ on the upper (lower) sides of the cut in $(z-1)^{\alpha-\frac{1}{2}}$.
 \\

For $x=\cos{\theta}\in(-1,1)$ with $\pi-\theta\ll1/\lvert\lambda\rvert^{1/3}$ and $\lvert\lambda\rvert\rightarrow\infty$,
\ba
\label{Dcut4}
{\mathsf D}_\lambda^\alpha(x) &\sim& \frac{\sqrt{\pi}}{\Gamma(\alpha)}\left(\frac{\lambda+\alpha}{2}\right)^{\alpha-1}\left(2(1+x)\right)^{-\alpha} \left(\frac{X''}{2}\right)^\frac{1}{2} \\ 
&& \times \left[-\sin{\pi\lambda}J_{\alpha-\frac{1}{2}}(X'')+\cos{\pi\lambda}Y_{\alpha-\frac{1}{2}}(X'')\right]\!, \nonumber \\
\label{Ccut4}
{\mathsf C}_\lambda^\alpha(x) &\sim&  \frac{\sqrt{\pi}}{\Gamma(\alpha)}\left(\frac{\lambda+\alpha}{2}\right)^{\alpha-1}\left(2(1+x)\right)^{-\alpha} \left(\frac{X''}{2}\right)^\frac{1}{2} \\
&& \times \left[\cos{\pi\lambda}J_{\alpha-\frac{1}{2}}(X'')+\sin{\pi\lambda}Y_{\alpha-\frac{1}{2}}(X'')\right]\!, \nonumber
\ea
where $X''=\sqrt{2(\lambda+\alpha)^2(1+x)}$, with uncertainties of relative order $1/\lvert\lambda\rvert^{2/3}$.
The results in Thms.\ 2 and 4 match for $1/\lvert\lambda\rvert\ll\lvert \sqrt{1+z}\rvert\ll 1/\lvert\lambda\rvert^{1/3}$. \\
\end{thm}

\noindent{\bf Derivation of Theorem 1:}

%\begin{proof}

Start with the following integral representation for $D_\lambda^\alpha(z)$ for $z\in{\mathbb C}$, $ \Re\lambda\geq 0$, and $ \Re(\lambda+2\alpha)>0$ \cite[Sec. 1\,(5)]{DFS}:
\be
\label{Dcontour}
D_\lambda^\alpha(x) = \frac{1}{2\pi i}e^{2\pi i\alpha}\int_{{\mathcal C}_+} dt\, t^{-\lambda-1}(t-z_+)^{-\alpha}(t-z_-)^{-\alpha},
\ee
where the integration contour ${\mathcal C}_+$ in the $t$ plane runs from $+\infty$, around the point $z_+$ in the positive sense, and back to $+\infty$, ${\mathcal C}_+=(\infty,z_++,\infty)$.  The factors $(t-z_\pm)^{-\alpha}$ are taken as cut in the $t$ plane from $z_\pm$ to $\infty$ along the directions defined by the lines from $t=0$ to $t=z_\pm$. The phases of the factors $t-z_\pm$ are defined as zero on the upper sides of the cuts for $\Im z>0$, and and elsewhere by continuation in $z$.

The function $C_\lambda^\alpha(z)$ has a similar integral representation for  $ \Re\lambda\geq 0$, and $ \Re(\lambda+2\alpha)>0$  \cite[Sec. 1\,(3)]{DFS}:
\be
\label{Ccontour}
C_\lambda^\alpha(z)=\frac{1}{2\pi i}e^{2\pi i\alpha}\int_{\mathcal C} dt\, t^{-\lambda-1}(t-z_+)^{-\alpha}(t-z_-)^{-\alpha}, 
\ee
where the contour ${\mathcal C}=(-\infty-i\epsilon,0+,-\infty+i\epsilon)$ runs around the negative $t$ axis in the positive sense.
In this representation the phases of the factors $(t-z_\pm)^{-\alpha}$ are defined separately for $\Im z\gtrless 0$, with, in both cases, the factors cut in the $t$ plane as above from $t=z_\pm$ to $\infty$, $0<{\rm arg}(t-z_\pm)<2\pi$. See \cite[Sec. 1\,(3)]{DFS} or \cite[Sec. 3.15.2\,(2)]{HTF}.

In these expressions, $z_\pm=z\pm\sqrt{z^2-1}$ with $\sqrt{z^2-1}$ cut in the $z$ plane from $z=1$ to $-\infty$, $-\pi< {\rm arg}z< \pi$. 
For $z$ in the upper (lower) half plane, $z_+$ is in the upper (lower) half plane outside the unit circle, while $z_-=1/z_+$ is  in the lower (upper) half plane inside the unit circle. For $z\in(-1,1)$, $z_\pm$ lie on the unit circle.
The singularities at $z_\pm$ pinch the contour ${\mathcal C}_+$ for $z\rightarrow \pm1$, so $D_\lambda^\alpha(z)$ has branch points at $\pm1$ and can be taken as cut from $\pm1$ to $-\infty$. Similarly, the singularities at $z_\pm$ pinch the contour ${\mathcal C}$ for $z\rightarrow -1$, so $C_\lambda^\alpha(z)$ has a branch point there and can be taken as cut from $-1$ to $-\infty$.

In treating the asymptotic properties of $C_\lambda^\alpha(z)$ and $D_\lambda^\alpha(z)$ in $\lambda$, we will take $\Re\alpha>0$ and $\Re \lambda\geq 0$. The integrands in Eqs.\ (\ref{Dcontour}) and (\ref{Ccontour}) are then singular at $t=0,\, z_+$, and $z_-$ and smaller in magnitude between, and vanish for $\lvert t\rvert\rightarrow\infty$, so there will be saddle points in the region of the singularities. If the contours ${\mathcal C}_+$ or $\mathcal{C}$ can be distorted to run through the saddle points  in the directions in which the integrands decrease most rapidly, the method of steepest descents provides an estimate of the integrals. This is valid provided the integrands are small  on the remainder of the contour and decrease rapidly for $\lvert t\rvert\rightarrow \infty$. 

To determine the location of the saddle points, write the integrands in Eqs.\ (\ref{Ccontour}) and (\ref{Dcontour}) as $e^{\Phi(t)}$, with
\be
\label{saddle_pts}
\Phi(t) = -(\lambda+1)\ln{t}-\alpha \ln{(t-z_+)}-\alpha \ln{(t-z_-)},
\ee
and require that $d\Phi/dt$ vanish as required for a stationary point. This gives the condition
\be
\label{saddle_pts2}
\frac{\lambda+1}{t}+\frac{\alpha}{t-z_+}+\frac{\alpha}{t-z_-} = 0.
\ee
For $\lvert\lambda\rvert$ large, the solutions must be close to $z_+$ or $z_-$. If those points are well separated, the solutions to order $1/\lvert\lambda \rvert$ are 
\be
\label{saddle_pts3}
t_\pm = z_\pm\left(1-\frac{\alpha}{\lambda+1}\right) +{\mathcal O}\left(\frac{\alpha^2}{(\lambda+1)^2}\right).
\ee
In general,
\be
\label{tplusminus_general}
t_\pm =\frac{ 1+\alpha'}{1+2\alpha'}\left[z\pm\sqrt{z^2-1+\left(\alpha'/(1+\alpha')\right)^2}\,\right],\quad \alpha'=\frac{\alpha}{\lambda+1},
\ee
so there are only the two saddle points $t_\pm$.

Next expand the exponent function $\Phi(t)$ in a Taylor series around the saddle points. To second order, for $\lvert\lambda\rvert$ large,
\be
\label{Phi_expanded}
\Phi(t) \approx \Phi(t_\pm)+\frac{\lambda^2}{\alpha z_\pm^2}(t-t_\pm)^2
\ee
near $t_\pm$. This gives the approximation
\ba
\label{saddle_integrand}
&& \frac{1}{2\pi i}e^{2\pi i\alpha} \left[z_+ \left(1-\frac{\alpha}{\lambda+1}\right)\right]^{-\lambda-1}\left(e^{i\pi}\frac{\alpha}{\lambda+1}z_+ \right)^{-\alpha} \\
&& \times\left(z_+-z_-\right)^{-\alpha}\int_{{\mathcal C}_+ } dt\, e^{\frac{1}{2}\frac{\lambda^2}{\alpha z_+ ^2}(t-t_+ )^2+\cdots} \nonumber
\ea
 for the integral in the neighborhood of the saddle point at $t_+$. The factors $e^{i\pi\alpha}\left(\alpha/(\lambda+1)\right)z_+$ and $(z_+-z_-)^{-\alpha}$ in this expression arise from the factors $(t_+-z_+)^{\alpha}$ and $(t_+-z_-)^{-\alpha}$ in the limit of large $\lvert\lambda \rvert$. A similar result holds  near $t_-$ with a different phase, $(t_--z_+)^{-\alpha}= e^{-i\pi\alpha}(z_+-z_-)^{-\alpha}[1+{\mathcal O}(1/\lvert\lambda\rvert]$.

The coefficient of $(t-t_+  )^2$ in the exponential  in the last factor in \eq{saddle_integrand} has phase $e^{2i\vartheta_+  }$, where
\be
\label{phaseTheta}
 \vartheta_+   =  {\rm arg}\,\lambda-\,{\rm arg}\,z_+  -\frac{1}{2}{\rm arg}\,\alpha
\ee
with  $-\pi/2\leq{\rm arg}\,\lambda\leq\pi/2$, $-\pi<{\rm arg}\,z_+<\pi$, and $-\pi/2<{\rm arg}\,\alpha<\pi/2$. 
With these ranges, the contour ${\mathcal C}_+ $ can be distorted to run through the saddle point in the direction with ${\rm arg}(t-t_+  ) =  \frac{\pi}{2}- \vartheta_+ $. The exponent is then real and negative, and the integration proceeds in the direction  of steepest descent away from the saddle.

The convergence of the integral away from the saddle point is rapid for $\lvert\lambda^2/\alpha z_+^2\rvert\gg 1$. Since the exact  integrand remains small on ${\mathcal C}_+ $ away from the saddle point, we can extend the integration on $t$ to $\pm\infty$ without changing the integral significantly. The result of the remaining Gaussian integral is just a factor $i\sqrt{2\pi\alpha z_+  ^2/\lambda^2}$, where the factor $e^{-i\vartheta_+ }$ from $dt$ has been absorbed.
Thus,  taking $\lvert\lambda\rvert$ large,
\ba
\label{D_integral}
D_\lambda^\alpha(z) &=& 2^{-\alpha}\left(e^\alpha \alpha^{-\alpha+1/2}/\sqrt{2\pi}\right){\Gamma(\alpha)}\lambda^{\alpha-1}(z^2-1)^{-\alpha/2}e^{i\pi\alpha} z_+^{-\lambda-\alpha} \\
&& \times  \left[1+{\mathcal O}(1/\lvert\lambda\rvert)\right]. \nonumber
\ea
The factor in parentheses is just Stirling's approximation for $1/\Gamma(\alpha)$, a known factor, in $D_\lambda^\alpha(z)$,  so 
\ba
\label{D_integral2}
D_\lambda^\alpha(z) &=& \frac{ 2^{-\alpha}}{\Gamma(\alpha)}\lambda^{\alpha-1}(z^2-1)^{-\alpha/2}e^{i\pi\alpha} z_+^{-\lambda-\alpha} \left[1+{\mathcal O}(1/\lvert\lambda\rvert)\right], \quad \lvert\lambda\rvert\rightarrow\infty, \\
&& \Re\lambda\geq 0,\quad \Re\alpha>0,\quad -\pi/2\leq {\rm arg}\lambda\leq\pi/2,\quad 0\leq{\rm arg}(z\pm 1)\leq\pi, \nonumber
\ea
in agreement with Eqs.\ (6.3) and (A5) in \cite{DFS}, but without the restriction on $\lambda$ noted there. This result holds in the complex $z$ plane cut from $z=1$ to $-\infty$.

 In the case of $C_\lambda^\alpha(z)$, we must distinguish the cases $\Im z>0$ and $\Im z<0$. For $\Im z>0$, the integral on the contour ${\mathcal C}_+$ reproduces the result for $D_\lambda^\alpha(z)$ in \eq{D_integral2}.  The integral on the ${\mathcal C}_-$ contour gives a similar result, with the replacement of $z_+$ by $z_-$ and an extra factor $e^{-i\pi\alpha}$ from the phase of the factor $(t_- -z_+)^{-\alpha}=e^{-i\pi\alpha}(z_+-z_-)^{-\alpha}\left(1+{\mathcal O}(1/\lambda)\right)$ in the integrand. 
 
 For $\Im z<0$, $(t_+-z_-)^{-\alpha} \rightarrow e^{-2\pi i \alpha}(z_+-z_-)^{-\alpha}$ for $\lvert\lambda\rvert$ large, and the factor $e^{i\pi\alpha}$ in \eq{D_integral2} from the ${\mathcal C}_+$ contour  is replaced by $e^{-i\pi\alpha}$. The contribution from ${\mathcal C}_-$ is unchanged.

 Combining the results for the ${\mathcal C}_+$ and ${\mathcal C}_-$ integrations, we find that for $\Im z\gtrless 0$, $ \lvert\lambda\rvert\rightarrow\infty$, $\Re\lambda\geq 0$, $\Re\alpha>0$, $ -\pi/2\leq {\rm arg}\lambda\leq\pi/2$,  and $0\leq{\rm arg}(z\pm 1)\leq\pi$, 

\ba
\label{Cresult1}
C_\lambda^\alpha(z) &=& \frac{2^{-\alpha}}{\Gamma(\alpha)}\lambda^{\alpha-1}(z^2-1)^{-\alpha/2}\left(e^{\pm i\pi\alpha} z_+^{-\lambda-\alpha} 
+  z_-^{-\lambda-\alpha} \right) 
\left[1+{\mathcal O}(1/\lvert\lambda\rvert)\right]  \\
\label{Cresult2}
&=&  \frac{2^{-\alpha}}{\Gamma(\alpha)}\lambda^{\alpha-1}(z^2-1)^{-\alpha/2}\left(e^{\pm i\pi\alpha} z_+^{-\lambda-\alpha} 
+  z_+^{\lambda+\alpha} \right)
 \left[1+{\mathcal O}(1/\lvert\lambda\rvert)\right].
\ea
%
%\end{proof}
This agrees with Eq.\ (A8) in \cite{DFS} and with Watson's result for $C_\lambda^\alpha(z)$, \cite[Sec. 2.3.2\,(17)]{HTF}. 

The earlier results for $C_\lambda^\alpha(x)$ and $D_\lambda^\alpha(x)$ were derived for $\lvert\lambda\rvert\rightarrow \infty$ along rays in the right-half $t$ plane with $\lvert\Im\lambda\rvert\rightarrow \infty$, $0<\lvert{\rm arg}\lambda\rvert<\pi/2$. The restrictions are not necessary, and the results continue to hold for ${\rm arg}\lambda=0$ and $\lvert{\rm arg}\lambda\rvert=\pi/2$.

The result for $C_\lambda^\alpha(z)$ must be interpreted with care.  Since $\lvert  z_+\rvert>1$ and $\lvert z_-\rvert=1/\lvert z_+\rvert<1$ for $z\not\in (-1,1)$, one of the two terms in \eq{Cresult2} usually becomes exponentially small relative to the other for $\lvert\lambda\rvert\rightarrow\infty$ and should be dropped relative to the uncertainties of $O(1/\lvert\lambda\rvert)$ in the dominant term. Thus, for   $\Re\lambda\rightarrow\infty$ with $\Im\lambda$ fixed, $\Re z>1$, and $\Im z\rightarrow 0$ from either above or below, the first, discontinuous, term in \eq{Cresult2} becomes exponentially small and should be dropped relative to the second. In fact, for $\lvert\Im{z}\rvert\rightarrow 0$, the saddle point at $z_+$ lies inside the contour for the $z_-$ integral, is inaccessible, and does not contribute to the final result. The asymptotic estimate for $C_\lambda^\alpha(z)$ is therefore continuous across the real axis as it should be.

This completes the derivation of the results in Theorem \ref{thm1D}. \\

%%%%%%%%%%%%%%%%%
%%%%%%%%%%%%%%%%%

\noindent {\bf Derivation of Theorem 2:} 

The functions ${\mathsf D}_\lambda^\alpha(x)$ and ${\mathsf C}_\lambda^\alpha(x)$ ``on the cut'', $x\in(-1,1)$, can be defined in terms of $D_\lambda^\alpha(z)$ for $z$ complex by \cite{Askey}
\ba
\label{Dcut}
{\mathsf D}_\lambda^\alpha(x) &=& -ie^{-i\pi\alpha}\left(e^{i\pi\alpha}D_\lambda^\alpha(x+i0)-e^{-i\pi \alpha}D_\lambda^\alpha(x-i0)\right), \\
\label{Ccut}
{\mathsf C}_\lambda^\alpha(x) &=&e^{-i\pi\alpha}\left( e^{i\pi\alpha}D_\lambda^\alpha(x+i0)+e^{-\pi i\alpha}D_\lambda^\alpha(x-i0)\right) \\
&=& C_\lambda^\alpha(x\pm i0).
\ea

For $z\in(-1,1)$,  take $z=\cos{\theta}$, $0<\theta<\pi$. Then $z_\pm=e^{\pm i\theta}$ and $\sqrt{z^2-1}=e^{\pm i\pi/2}\sin{\theta}$ for $\Im z\gtrless 0$, and Eqs.\ (\ref{D_integral2}) and (\ref{Cresult2}) give
\ba
{\mathsf D}_\lambda^\alpha(\cos{\theta})  &=& i\frac{2^{-\alpha}}{\Gamma(\alpha)}\lambda^{\alpha-1}(\sin{\theta})^{-\alpha}\left(-e^{-i(\lambda+\alpha)\theta+i\pi\alpha/2}+e^{i(\lambda+\alpha)\theta-i\pi\alpha/2}\right) \\
&& \times\left[1+{\mathcal O}(1/\lvert\lambda\rvert)\right]  \nonumber \\
\label{Dcut_result}
&=&-\frac{2^{-\alpha+1}}{\Gamma(\alpha)}\lambda^{\alpha-1}(\sin{\theta})^{-\alpha}\sin{((\lambda+\alpha)\theta-\pi\alpha/2)}\left[1+{\mathcal O}(1/\lvert\lambda\rvert)\right], \nonumber \\
{\mathsf C}_\lambda^\alpha(\cos{\theta)} &=& \frac{2^{-\alpha}}{\Gamma(\alpha)}\lambda^{\alpha-1}(\sin{\theta})^{-\alpha}\left(e^{i\pi\alpha/2} e^{-i(\lambda+\alpha)\theta}  + e^{-i\pi\alpha/2}  e^{i(\lambda+\alpha)\theta} \right) \\
&& \times \left[1+{\mathcal O}(1/\lvert\lambda\rvert)\right] \nonumber \\
\label{Ccut_result}
&=& \frac{2^{-\alpha+1}}{\Gamma(\alpha)}\lambda^{\alpha-1}(\sin{\theta})^{-\alpha}\cos{\left((\lambda+\alpha)\theta-\pi\alpha/2\right)}\left[1+{\mathcal O}(1/\lvert\lambda\rvert)\right]. \nonumber
\ea

A question now is how large $|\lambda|$ must actually be for this behavior to hold. It follows from the expressions for $t_\pm$ in \eq{tplusminus_general} that the  asymptotic limit for the saddle points that used in the calculations requires that $\sqrt{z^2-1}\gg \alpha'\approx\alpha/\lambda$. Furthermore, the points $t_\pm$ or $z_\pm$ must be separated widely enough that the integration over one saddle is not influenced by the presence of the second. 

The convergence of the saddle point integrals is determined by the coefficient in the exponential in the integral in \eq{saddle_integrand}. Convergence on the right scale requires that the distance between the points be much larger than the sum of the distances over which the saddle point integrations converge, given by the scale factors $\sqrt{\lvert 2\alpha z_\pm^2/\lambda^2\rvert}$ in the Gaussian integrands. This gives the condition
\be
\label{theta_cond1}
|z_+-z_-| = 2\sin{\theta} \gg \sqrt{|2\alpha z_+^2/\lambda^2}| + \sqrt{|2\alpha z_-^2/\lambda^2}| = 2\sqrt{|2\alpha/\lambda^2|}
\ee
for $z\in(-1,1)$, so requires that
\be
\label{lambda_limit}
\lvert\lambda\rvert \gg \left|\sqrt{\alpha}\big/\sqrt{z^2-1}\right|.
\ee
This is the same as the condition used in the derivation of $t_\pm$ given above up to a factor $\sqrt{\alpha}$. 

 For fixed large $|\lambda|$,  \eq{lambda_limit} bounds $|z^2-1|$ away from 1. The saddle points $t_\pm$ merge for $\theta\rightarrow 0$ ($z\rightarrow 1$) and cannot be treated as independent in the steepest-descent calculations which lead to the results above. For $z\rightarrow1$, the saddle points coalesce into a single saddle between $t=0$ and $t=1$, and an integration as in \eq{saddle_integrand} with $\alpha\rightarrow 2\alpha$ reproduces the correct asymptotic limit  $C_\lambda^\alpha(1)\sim \lambda^{2\alpha-1}/\Gamma(2\alpha)$. For $z\rightarrow -1$, the points $z_\pm$ pinch the contour $\mathcal{C}$, and the result is singular, $C_\lambda^\alpha(z)\propto \left((z+1)/2\right)^{-\alpha+1/2}$. \\
 
 %%%%%%%%%%%%%%%%%
 %%%%%%%%%%%%%%%%%

\noindent {\bf Derivation of Theorem 3:}

To treat the limit $z\rightarrow 1$, we use a different technique developed in \cite[Secs. IIA, IIB]{LDlegendre}. We start with the standard hypergeometric expression for $C_\lambda^\alpha(z)$ written in a more useful form,
\ba
\label{C_2F1}
C_\lambda^\alpha(z) &=& \frac{\Gamma(\lambda+2\alpha)}{\Gamma(\lambda+1)\Gamma(2\alpha)}\, _2F_1\left(-\lambda,\lambda+2\alpha;\alpha+\frac{1}{2};\frac{1-z}{2}\right) \\
\label{alt_C}
&=& \frac{\Gamma(\lambda+2\alpha)}{\Gamma(\lambda+1)\Gamma(2\alpha)}\left(\frac{1+z}{2}\right)^{-\alpha+\frac{1}{2}}\\ && \times \, _2F_1\left(\lambda+\alpha+\frac{1}{2},-\lambda-\alpha+\frac{1}{2};\alpha+\frac{1}{2};\frac{1-z}{2}\right). \nonumber
\ea
We next introduce the Barnes-type representation   \cite[Sec. 2.3.3\,(15)]{HTF} for the type of hypergeometric function  that appears in \eq{alt_C} and will be encountered again later,
\ba
&& _2F_1\left(b+\frac{1}{2},-b+\frac{1}{2};\nu+1;\frac{u}{2}\right)  \\ 
\label{2F1_Barnes}
&& =  \Gamma(\nu+1)\,\frac{1}{2\pi i}\int_{{\mathcal C}_B} ds\frac{\Gamma(b+\frac{1}{2}+s)\Gamma(-b+\frac{1}{2}+s)}{\Gamma(b+\frac{1}{2})\Gamma(-b+\frac{1}{2})} \frac{\Gamma(-s)}{\Gamma(\nu+1+s)}\left(-\frac{u}{2}\right)^s.  \nonumber
\ea
The contour ${\mathcal C}_B$ in the Barnes' representation initially runs from $-i\infty$ to $+i\infty$  in the $s$ plane, staying to the right of the poles of the factors $\Gamma(b+\frac{1}{2}+s)$ and $\Gamma(-b+\frac{1}{2}+s)$ in the integrand, and to the left of the poles of $\Gamma(-s)$, but it can be deformed to run around the positive real axis, ${\mathcal C}_B=(\infty,0-,\infty)$ with the same restrictions.

Expanding   the ratios of $b$-dependent gamma functions in the first line of \eq{2F1_Barnes} in inverse powers of $b$, assumed large, using Stirling's approximation for the gamma function, and writing the powers of $s$ that appear in terms of combinations of the form $1\cdot s(s-1)\cdots(s-k)$ gives a series
\ba
&& _2F_1\left(b+\frac{1}{2},-b+\frac{1}{2};\nu+1;\frac{u}{2}\right)  \label{C_Barnes_series} \\
&& \qquad  =  \Gamma(\nu+1) \frac{1}{2\pi i}\int_{{\mathcal C}_B}ds\frac{\Gamma(-s)}{\Gamma(\nu+1+s)}\left(\frac{U^2}{4}\right)^s  \nonumber \\
&& \qquad  \times\left\{1-\frac{1}{b^2}\left(\frac{s}{4}+s(s-1)+\frac{1}{3}s(s-1)(s-2)\right)  
+{\mathcal O}\left(\frac{1}{b^4}\right)\right\}, \nonumber
\ea
where $U=\sqrt{2b^2u}$.

With the choice of the deformed contour ${\mathcal C}_B$ above, the integrals that remain in \eq{2F1_Barnes_expanded},  are expressible in terms of Bessel functions through a Barnes' representation for the latter which uses the same contour ${\mathcal C}_B$,
\be
\label{Bessel_Barnes}
\left(\frac{U}{2}\right)^{-\nu}\,J_\nu(U)=\frac{1}{2\pi i}\int_{{\mathcal C}_B} ds\, \frac{\Gamma(-s)}{\Gamma(\nu+s+1)}\left(\frac{U^2}{4}\right)^s.
\ee
The first term in the series in \eq{2F1_Barnes_expanded} gives $(U/2)^{-\nu}J_\nu(U)$. After combining the factors $s(s-1)\cdots(s-k)$, $k=0,\,1,\cdots$, with $\Gamma(-s)$ to get $(-1)^{k+1}\Gamma(-s+k+1)$, we can shift the contour of integration to the right to run just to the left of the pole at $s=k+1$. The replacement of $s$ by $s'=s-k-1$ then gives
\ba
\label{J_shifted}
&&\frac{1}{2\pi i}\int_{{\mathcal C}_B}ds\,s(s-1)\cdots(s-k)\frac{\Gamma(-s)}{\Gamma(s+\nu+1)}\left(\frac{U}{2}\right)^s  \\ 
&&\quad = \frac{(-1)^{k+1}}{2\pi i}\int_{{\mathcal C}'_B}ds'\,\frac{\Gamma(-s')}{\Gamma(s'+\nu+k+2)} \left(\frac{U^2}{4}\right)^{s'+k+1} \nonumber \\
&& \quad 
=(-1)^{k+1}\left(\frac{U}{2}\right)^{-\nu+k+1 }J_{\nu+k+1}(U) \nonumber
\ea
for the following terms, with $k=0,\,1,\cdots$.

The use of Stirling's approximation, itself only an asymptotic expansion, is not justified on the entire integration contour, and the result from \eq{C_Barnes_series} gives only an asymptotic series for the hypergeometric function,
\ba
\label{2F1Bessel}
&& _2F_1\left(b+\frac{1}{2},-b+\frac{1}{2};\nu+1;\frac{u}{2}\right) = \Gamma(\nu+1)\left(\frac{U}{2}\right)^{-\nu}  \\ 
&& \times \left\{J_\nu(U)+\frac{1}{b^2}\left[\frac{1}{4}\frac{U}{2}J_{\nu+1}(U) 
- \left(\frac{U}{2}\right)^2J_{\nu+2}(U)+\frac{1}{3}\left(\frac{U}{2}\right)^3J_{\nu+3}(U)\right] \right. \nonumber \\
&& + \left. {\mathcal O}\left(\frac{1}{b^4}\right) \right\}. \nonumber
\ea

The use of this expression in \eq{C_Barnes_series} with $b=\lambda+\alpha$ and $\nu=\alpha-\frac{1}{2}$  gives an asymptotic series for $C_\lambda^\alpha(z)$ in powers of $1/(\lambda+\alpha)^2$. With $Z=\sqrt{2(\lambda+\alpha)^2(1-z)}$,
\ba
\label{Bessel_rep}
C_\lambda^\alpha(z) &=&  \frac{\Gamma(\lambda+2\alpha)\Gamma(\alpha+\frac{1}{2})}{\Gamma(\lambda+1)\Gamma(2\alpha)}\left(\frac{1+z}{2}\right)^{-\alpha+\frac{1}{2}}\left(\frac{Z}{2}\right)^{-\alpha+\frac{1}{2}}\Bigg\{ J_{\alpha-\frac{1}{2}}(Z)   \\
&&  \!\!\!\!\! \!\!\!\!\!\!  +\left. \frac{1}{(\lambda+\alpha)^2}\left[\frac{1}{4}\frac{Z}{2}J_{\alpha+\frac{1}{2}}(Z)-\left(\frac{Z}{2}\right)^2J_{\alpha+\frac{3}{2}}(Z)+\frac{1}{3}\left(\frac{Z}{2}\right)^3J_{\alpha+\frac{5}{2}}(Z)\right] \right. \nonumber \\
&& \!\!\!\!\!\!\!\!\!\!\!  + \left. {\mathcal O}\left(\frac{1}{(\lambda+\alpha)^4}\right)\right\} \nonumber
\ea
for $\lambda+\alpha$ large and $Z$ fixed. This series is useful more generally for $|\lambda+\alpha |\gg 1$ and $\lvert 1-z\rvert\ll 1$.

We can obtain a closely-related series using the same technique starting with \eq{C_2F1} and expanding in terms of the parameter $\lambda(\lambda+2\alpha)=(\lambda+\alpha)^2-\frac{1}{4}$. This approach was used in \cite[Sec. IIA]{LDlegendre} in our treatment of Bessel-function expansions for the associated Legendre functions $P_j^{-\mu}(z)$.  The result is
\ba
C_\lambda^\alpha(z) &=& \frac{\Gamma(\lambda+2\alpha)\Gamma(\alpha+\frac{1}{2})}{\Gamma(\lambda+1)\Gamma(2\alpha)}\left(\frac{Y}{2}\right)^{-\alpha+\frac{1}{2}}\left\{J_{\alpha-\frac{1}{2}}(Y) \right.  
\label{Bessel_rep2} \\
&&  \left. +\frac{1}{\lambda(\lambda+2\alpha)}\left[-\frac{2\alpha+1}{2}\left(\frac{Y}{2}\right)^2J_{\alpha+\frac{3}{2}}(Y)+\frac{1}{3}\left(\frac{Y}{2}\right)^3J_{\alpha+\frac{5}{2}}(Y)\right]  \right. \nonumber \\
&& \left. +{\mathcal O}\left(\frac{1}{(\lambda(\lambda+2\alpha)^2}\right)\right\}, \nonumber
\ea
where $Y=\sqrt{2\lambda(\lambda+2\alpha)(1-z)}$. 

The series in \eq{Bessel_rep}, here obtained directly, is equivalent to that obtained by expanding the powers of $\lambda(\lambda+2\alpha)$ in the coefficients and the argument of the Bessel functions in \eq{Bessel_rep2} in terms of the simpler variable $(\lambda+\alpha)^2$, and the prefactor $\left((1+z)/2\right)^{-\alpha+\frac{1}{2}}$ in powers of $(1-z)$. The difference in the leading terms is unimportant for $\lvert\lambda\rvert\gg1$ and $\lvert1-z\rvert\ll1$.

To connect this result to the asymptotic expression for $C_\lambda^\alpha(z)$ in \eq{Ccut_result}  for $\lvert\lambda\rvert\gg 1$, we consider the case in which only the leading term in the asymptotic series in  \eq{Bessel_rep} is important. The result in \eq{Ccut_result} is valid for $\sqrt{\lvert 1-z\rvert}\gg 1/\lvert\lambda\rvert$, which requires that $Z\gg 1$. Despite the appearance of powers of $Z$ in the correction terms, this is allowed provided that the  corrections to the leading term are small. The Bessel functions are all of the same general magnitude for $Z$ large, so the term in $Z^3$ in the second term in the series is dominant for $Z\gg 1$, and the condition for the the $\lvert1/(\lambda+\alpha)^2\rvert$ correction to the leading term to be small is
\be
\label{Z_condition}
\left\lvert\frac{1}{(\lambda+\alpha)^2}\left(\frac{Z}{2}\right)^3\right\rvert =\vert\lambda+\alpha\rvert \left\lvert\frac{1-z}{2}\right\rvert^{3/2} \ll 1.
\ee
Under this condition, the following terms in the series in \eq{Bessel_rep} are also initially small.

For $z=\cos{\theta}$ on $-1<z<1$, this requires that $\theta$ be small, with $\theta\ll 2/(\lambda+\alpha)^{\frac{1}{3}}$. 
 In the limit of large $Z$, Hankel's expansion for the Bessel functions  \cite[Sec. 10.17(i)]{DLMF}  gives     
\ba
J_\nu(Z) &=& \sqrt{\frac{2}{\pi Z}}\left\{\cos\left(Z-\left(\nu+\frac{1}{2}\right)\frac{\pi}{2}\right)\left[1+{\mathcal O}\left(\frac{1}
{Z^2}\right)\right]  \right.  \label{Bessel_cos} \\
&&  \left.  -\sin\left(Z-\left(\nu+\frac{1}{2}\right)\frac{\pi}{2}\right)\left[\frac{4\nu^2-1}{8Z}+{\mathcal O}\left(\frac{1}{Z^3}\right)\right] \right\}. \nonumber
\ea
The leading term in \eq{Bessel_rep} therefore has the asymptotic limit
\ba
\label{C_Bessel_asym}
C_\lambda^\alpha(z) &\sim&  \frac{\Gamma(\lambda+2\alpha)\Gamma(\alpha+\frac{1}{2})}{\Gamma(\lambda+1)\Gamma(2\alpha)}\left(\frac{1+z}{2}\right)^{-\alpha+\frac{1}{2}}\left(\frac{Z}{2}\right)^{-\alpha+\frac{1}{2}} \\
&& \times \sqrt{\frac{2}{\pi Z}} \cos\left(Z-\frac{\pi\alpha}{2}\right). \nonumber
\ea
The corrections are of relative order $1/Z\propto \lvert\lambda\rvert^{-\frac{2}{3}}$.

Expanding the leading ratio of gamma functions in terms of $\lambda+\alpha$, assumed large, and expressing the result in terms of $\theta\ll 1/\lvert\lambda\rvert^{1/3}$, we obtain
\ba
\label{C_theta_form}
C_\lambda^\alpha(\cos{\theta}) &=& \left(\frac{\lambda+\alpha}{2}\right)^{\alpha-1}\frac{1}{\Gamma(\alpha)}\left(\sin{\theta}\right)^{-\alpha}\cos\left((\lambda+\alpha)\theta-\frac{\pi\alpha}{2}\right) \\
&& \times \left[1+{\mathcal O}\left(1/\lvert\lambda\rvert^{2/3}\right)\right]. \nonumber
\ea
This is equivalent to the expression in \eq{Ccut_result}  for $\lvert\lambda\rvert\gg 1$  as required for the the correction terms in \eq{Bessel_rep} to be negligible, so the two expressions connect smoothly in their overlapping region of validity, $1/\lvert\lambda\rvert\ll\theta\ll1/\lvert\lambda\rvert^{1/3}$.

To obtain an asymptotic Bessel-function series for $D_\lambda^\alpha(z)$ for $z \approx 1$, we use the relation \cite[Sec. 3\,(5)]{DFS}
\ba
\label{DFS3.5}
D_\lambda^\alpha(z) &=& \frac{1}{2}e^{i\pi\alpha}\frac{1}{\cos{\pi\alpha}} \\
&\times& \left\{C_\lambda^\alpha(z)-2^{-2\alpha+1}(z^2-1)^{\alpha+\frac{1}{2}}\frac{\Gamma(-\alpha+1)\Gamma(\lambda+2\alpha)}{\Gamma(\alpha)\Gamma(\lambda+1)}C_{\lambda+2\alpha-1}^{-\alpha+1}(z)\right\} \nonumber
\ea
to express $D_\lambda^\alpha(z)$ in terms of Gegenbauer functions of the first kind. The function $C_\lambda^\alpha(z)$ can be approximated using the series in \eq{Bessel_rep}. The modified indices $\lambda'=+2\alpha-1$ and $\alpha'=-\alpha+1$ in the second Gegebauer function give $\lambda'+\alpha'=\lambda+\alpha$ so we may use the same series for this function, but with the index $\alpha-\frac{1}{2}$ on the Bessel functions and their coefficients replaced by $\alpha'-\frac{1}{2}=-\alpha+\frac{1}{2}$.

 We begin with $z\in(1,\infty)$, where $Z\rightarrow iZ'$, $Z'=\sqrt{2(\lambda+\alpha)^2(z-1)}$ and $Z^{\mp\nu}J_{\pm\nu+n}(Z)\rightarrow (-1)^n(Z')^{\mp\nu}I_{\pm\nu+n}(Z')$, with $I_\mu(z)$  the modified or hyperbolic Bessel function of the first kind. This gives as the leading term 
\ba
\label{D_Bessel1}
D_\lambda^\alpha(z) &=& \frac{\sqrt{\pi}}{2}e^{i\pi\alpha}\frac{1}{\Gamma(\alpha)}\frac{1}{\cos{\pi\alpha}} \\ 
&& \times \left\{ \frac{\Gamma(\lambda+2\alpha)}{\Gamma(\lambda+1)}2^{-2\alpha+1}\left(\frac{z+1}{2}\right)^{-\alpha+\frac{1}{2}}\left(\frac{Z'}{2}\right)^{-\alpha+\frac{1}{2}}I_{\alpha-\frac{1}{2}}(Z')\right.  \nonumber \\
&& \left. -(z^2-1)^{-\alpha+\frac{1}{2}}\left(\frac{z+1}{2}\right)^{\alpha-\frac{1}{2}}\left(\frac{Z'}{2}\right)^{\alpha-\frac{1}{2}}I_{-\alpha+\frac{1}{2}}(Z') \right. \nonumber \\
&& \left. +{\mathcal O}\left(\frac{1}{(\lambda+\alpha)^2}\right)\right\}. \nonumber
\ea
The higher-order terms in the series are negligible for $\lvert z-1\rvert\ll\vert\lambda\rvert^{-\frac{2}{3}}$ for $\lvert\lambda\rvert\rightarrow\infty$.

The function $I_{-\alpha+\frac{1}{2}}(Z')$ can be eliminated in terms of the Macdonald function $K_{\alpha-\frac{1}{2}}(Z')$ through the relation \cite[Sec. 10.27.4)]{DLMF}
\be
\label{Kfunction}
K_\nu(z) = \frac{\pi}{2}\frac{1}{\sin{\pi\nu}}\big[I_{-\nu}(z)-I_\nu(z)\big]
\ee
with $\nu=\alpha-\frac{1}{2}$. After some rearrangement of the coefficients, this gives
\ba
D_\lambda^\alpha(z) &\sim& \frac{1}{\sqrt{\pi}}e^{i\pi\alpha}\frac{1}{\Gamma(\alpha)}2^{-\nu}(\lambda+\alpha)^\nu (z^2-1)^{-\nu/2}\left(\frac{z+1}{2}\right)^{\nu/2}\bigg\{ K_\nu(Z')   \\
\label{D_Bessel2}
&&\!\!\!\!\!\!\!\!\!\!\!\!\!\!\! + \left.\frac{\pi}{2}\frac{1}{\sin{\pi\nu}}\left[1-
\frac{\Gamma(\lambda+2\alpha)}{\Gamma(\lambda+1)} (\lambda+\alpha)^{-2\alpha+1}\left(\frac{z+1}{2}\right)^{-\nu}\right]I_\nu(Z')+\cdots\right\}. \nonumber
\ea
Since $\Gamma(\lambda+2\alpha)/\Gamma(\lambda+1)= (\lambda+\alpha)^{2\alpha-1}[1+{\mathcal O}(1/(\lambda+\alpha)^2]$ while $\big(z+1)/2\big)^{-\nu}\sim 1-\nu(z-1)/2+\cdots$, the coefficient of $I_\nu(Z')$ in this expression vanishes up to terms of order $1/(\lambda+\alpha)^2$  and $(z-1)/2\ll(\lambda+\alpha)^{-\frac{2}{3}}$ over its range of validity. The overall factor $\big((z+1)/2\big)^{-\nu/2}$ can also be dropped to leading order, and
\be
\label{D_Bessel3}
D_\lambda^\alpha(z) \sim \frac{1}{\sqrt{\pi}}e^{i\pi\alpha}\frac{1}{\Gamma(\alpha)}2^{-\nu}(\lambda+\alpha)^\nu (z^2-1)^{-\nu/2} K_\nu(Z'),\quad \nu=\alpha-\frac{1}{2}.
\ee

We obtain the asymptotic forms of the Gegenbauer functions on the cut for $x=\cos{\theta}\approx 1$ using  Eqs.\ (\ref{Dcut}) and (\ref{Ccut}) and the relations
\be
\label{K(x+-i0)}
K_\nu(e^{\pm i\pi/2}z) = \mp\frac{i\pi}{2}e^{\mp i\pi\nu/2}\left[J_\nu(z)\mp iY_\nu(z)\right],
\ee
This gives
\ba
\label{Dcut3}
{\mathsf D}_\lambda^\alpha(x) &\sim& -\sqrt{\pi}\frac{1}{\Gamma(\alpha)}2^{-\nu}(\lambda+\alpha)^{\nu}(1-x^2)^{-\nu/2}Y_\nu(Z)+\cdots, \\
\label{Ccut3}
{\mathsf C}_\lambda^\alpha(x) &\sim&\sqrt{\pi} \frac{1}{\Gamma(\alpha)}2^{-\nu}(\lambda+\alpha)^{\nu}(1-x^2)^{-\nu/2}J_\nu(Z)+\cdots, 
\ea
$\nu=\alpha-\frac{1}{2}$, where the uncertainties in these expressions in their range of validity are order $1/\lambda^{2/3}$  for $\lvert\lambda\rvert\rightarrow\infty$.

Using Hankel's expansions of $J_\nu(Z)$ and $Y_\nu(Z)$ for $Z$ large \cite[Sec. 10.17(i)]{DLMF} and expressing the results in terms of $\theta$, with $x=\cos{\theta}$, $Z=(\lambda+\alpha)\sqrt{2(1-\cos{\theta}}=(\lambda+\alpha)\theta[1-\theta^2/24+\cdots]$, and $\theta\ll \lvert\lambda\rvert^{-\frac{1}{3}}$, these relations give
\ba
\label{Dcut4}
\ \ {\mathsf D}_\lambda^\alpha(\cos{\theta}) &\sim& -\frac{1}{\Gamma(\alpha)}2^{-\alpha+1}(\lambda+\alpha)^{\alpha-1}(\sin{\theta})^{-\alpha}\sin\left((\lambda+\alpha)\theta-\frac{\alpha\pi}{2}\right)+\cdots, \\
\label{Ccut4}
{\mathsf C}_\lambda^\alpha(\cos{\theta}) &\sim& \frac{1}{\Gamma(\alpha)}2^{-\alpha+1}(\lambda+\alpha)^{\alpha-1}(\sin{\theta})^{-\alpha}\cos\left((\lambda+\alpha)\theta-\frac{\alpha\pi}{2}\right)+\cdots,
\ea
to leading order in $\lvert\lambda\rvert$, in agreement with the results in Eqs.\ (\ref{Dcut_result}) and (\ref{Ccut_result}) for $1/\lvert\lambda\rvert\ll\theta\ll1/\lvert\lambda\rvert^{1/3}$. \\

%%%%%%%%%%%%%%
%%%%%%%%%%%%%%

\noindent{\bf Derivation of Theorem 4:}

Begin with the following hypergeometric representaion of $D_\lambda^\alpha(z)$ for $\lvert z\rvert$ large\cite{DFS,Askey},
\ba
\label{D_hyper1}
D_\lambda^\alpha(z) &=& e^{i\pi\alpha}\left(2(z-1)\right)^{-\lambda-2\alpha}\frac{\Gamma(\lambda+2\alpha)}{\Gamma(\lambda+\alpha+1)\Gamma(\alpha)} \\
&& \times  _2F_1\left(\lambda+2\alpha,\lambda+\alpha+\frac{1}{2};2\lambda+2\alpha+1;\frac{2}{1-z}\right). \nonumber
\ea
This can be converted using standard linear transformations \cite[Sec. 15.8]{DLMF} to a form useful for complex $z$ near $-1$,
\ba
 \label{D_hyper2}
D_\lambda^\alpha(z) &=& \frac{1}{\sqrt{\pi}}e^{i\pi\alpha}2^{-2\alpha}\frac{1}{\Gamma(\alpha)}e^{\mp i\pi(\lambda+2\alpha)} \left[ \frac{\Gamma(\lambda+2\alpha)\Gamma(-\alpha+\frac{1}{2})}{\Gamma(\lambda+1)} \right. \\&& \times \left.\left(\frac{1-z}{2}\right)^{-\alpha+\frac{1}{2}} \, _2F_1\left(-\lambda-\alpha+\frac{1}{2},\lambda+\alpha+\frac{1}{2};\alpha+\frac{1}{2};\frac{1+z}{2}\right) \right.
 \nonumber \\
&+& \left. e^{\pm i\pi(\alpha-\frac{1}{2})}\Gamma\left(\alpha-\frac{1}{2}\right)\left(\frac{1+z}{2}\right)^{-\alpha+\frac{1}{2}} \right. \nonumber \\
&& \times \left. _2F_1\left(-\lambda-\alpha+\frac{1}{2},\lambda+\alpha+\frac{1}{2};-\alpha+\frac{3}{2};\frac{1+z}{2}\right)  \right],  \nonumber
\ea
where the + and - signs hold for $z$ on the upper (lower) sides of the cut in $(z-1)^{\alpha-\frac{1}{2}}$.

Upon using the asymptotic Bessel-function approximation in \eq{2F1Bessel} for the hypergeometric functions in leading order and expanding the ratio of gamma functions in the first term for $\lvert\lambda+\alpha\rvert\gg1$, this reduces in leading order to
\ba
\label{DBessel(Z'')}
D_\lambda^\alpha(z) &\sim& e^{i\pi\alpha}2^{-\alpha}(\lambda+\alpha)^{\alpha-1}\frac{\sqrt{\pi}}{\sin{\pi(\alpha-\frac{1}{2})}}\frac{1}{\Gamma(\alpha)}\left(2(1+z)\right)^{-\alpha/2}\left(\frac{Z''}{2}\right)^{\frac{1}{2}} \\
&& \!\!\!\!\!\!\!\!\!\! \!\!\! \times e^{\mp i\pi(\lambda+2\alpha)}\left\{-\left(\frac{1-z}{2}\right)^{-\alpha+\frac{1}{2}}J_{\alpha-\frac{1}{2}}(Z'')+e^{\pm i\pi(\alpha-\frac{1}{2})}J_{-\alpha+\frac{1}{2}}(Z'')\right\}, \nonumber
\ea
where $Z''=\sqrt{2(\lambda+\alpha)^2(1+z)}$. The factor $\left((1-z)/2\right)^{-\alpha+\frac{1}{2}}$ differs from 1 only by corrections of order $1/\lvert\lambda\rvert^\frac{2}{3}$ in the region in which the leading-order approximation is valid, so can be replaced by 1 for $\lvert\lambda\rvert\gg 1$.

For $z=x\in(-1,1)$ real and close to $-1$, with $x=\cos{\theta}$, $\theta\approx\pi$,  the relations in Eqs.\  (\ref{Dcut}),  (\ref{Ccut}), and (\ref{DBessel(Z'')}), give the asymptotic forms of the Gegenbauer functions ${\mathsf D}_\lambda^\alpha$ and ${\mathsf C}_\lambda^\alpha(x)$ ``on the cut'" for $\lvert\lambda\rvert\gg 1$.  Calculating the discontinuities specified in the first two equations and replacing $J_{-\alpha+\frac{1}{2}}(x)$ by the Bessel function of the second kind,
\be
\label{Ydefined}
Y_{\alpha+\frac{1}{2}}(x) = \frac{1}{\sin{\pi(\alpha-\frac{1}{2}}}\left(J_{\alpha-\frac{1}{2}}(x)\cos{\left(\pi(\alpha-\frac{1}{2}\right)} - J_{-\alpha+\frac{1}{2}}(x)\right)
\ee
gives the relations in Thm. 4,
\ba
\label{DcutY}
{\mathsf D}_\lambda^\alpha(x)  &\sim& \frac{\sqrt{\pi}}{\Gamma(\alpha)}\left(\frac{\lambda+\alpha}{2}\right)^{\alpha-1}\left(2(1+x)\right)^{-\alpha/2}\left(\frac{X''}{2}\right)^\frac{1}{2}  \\
&& \times \left[-\sin{\pi\lambda}J_{\alpha-\frac{1}{2}}(X'')+\cos{\pi\lambda}Y_{\alpha-\frac{1}{2}}(X'')\right]\!, \nonumber \\
\label{CcutY}
 {\mathsf C}_\lambda^\alpha(x)  &\sim&  \frac{\sqrt{\pi}}{\Gamma(\alpha)}\left(\frac{\lambda+\alpha}{2}\right)^{\alpha-1}\left(2(1+x)\right)^{-\alpha/2}\left(\frac{X''}{2}\right)^\frac{1}{2} \\
&& \times \left[\cos{\pi\lambda}J_{\alpha-\frac{1}{2}}(X'')+\sin{\pi\lambda}Y_{\alpha-\frac{1}{2}}(X'')\right], \nonumber
\ea
with $X''=\sqrt{2(\lambda+\alpha)^2(1+x)}\approx (\lambda+\alpha)(\pi-\theta)$. 

For $1/\lvert\lambda\rvert\ll\pi-\theta\ll 1/\lvert\lambda\rvert$, the results in  Eqs.\ (\ref{DcutY}) and (\ref{thm2D}) and in (\ref{CcutY}) and   (\ref{thm2C}) are in their common ranges of validity and should agree. $X''$ is large in this region, and the agreement is easily shown using Hankel's asymptotic expressions for the Bessel functions \cite[Sec. 10.17(i)]{DLMF} and noting that $\sin{\theta}=\sin{(\pi-\theta)}\approx \pi-\theta$ in this region. \\

 %%%%%%%%%%%%%%%%%
 %%%%%%%%%%%%%%%%%

\noindent{\bf Remarks:}

 In their discussion of the asymptotics of  the associated Legendre functions
 \be
 \label{PfromC}
  P_\nu^{-\mu}(z)= \frac{2^\mu}{\sqrt{\pi}}\frac{\Gamma(\mu+\frac{1}{2})\Gamma(\nu-\mu+1)}{\Gamma(\nu+\mu+1)}\left(z^2-1\right)^{\mu/2} C_{\nu-\mu}^{\mu+\frac{1}{2}}(z) \nonumber
  \ee
  and
  \be
 \label{QfromD}
  Q_\nu^{-\mu}(z)= 2^\mu\sqrt{\pi}e^{-2\pi i(\mu+\frac{1}{4})}\frac{\Gamma(\mu+\frac{1}{2})\Gamma(\nu-\mu+1)}{\Gamma(\nu+\mu+1)}\left(z^2-1\right)^\frac{\mu}{2}D_{\nu-\mu}^{\mu+\frac{1}{2}}(z) \nonumber
  \ee
  for $\lvert\nu\rvert\rightarrow\infty$,  Cohl, Dang, and Dunster \cite[Secs. 2.3.1, 2.4.1]{cohl} use uniform asymptotic expressions in terms of Bessel functions which hold quite generally  \cite[Sec.14.15\,(11-14)]{DLMF}. These involve arguments $(\mu+\frac{1}{2})\theta$ in the Bessel functions and pre-factors proportional to $\sqrt{\theta/\sin{\theta}}$ for $z=\cos{\theta}$ or $\sqrt{\theta/\sinh{\theta}}$ for $z=\cosh{\theta}$. For example,  the Ferrers functions ${\mathsf P}_\nu^{-\mu}(\cos{\theta})$ and ${\mathsf Q}_\nu^{-\mu}(\cos{\theta})$ have the asymptotic forms
 \ba
 \label{FerrersPlimit}
 {\mathsf P}_\nu^{-\mu}(\cos{\theta}) &=& \frac{1}{\nu^\mu}\left(\frac{\theta}{\sin{\theta}}\right)^\frac{1}{2} \left[J_\mu\left(\left(\nu+\frac{1}{2}\right)\theta\right)+{\mathcal O}\left(\frac{1}{\nu}\right){\rm env}J_\mu\left(\left(\nu+\frac{1}{2}\right)\theta\right)\right], \nonumber \\
 \label{FerrersQlimit}
 {\mathsf Q}_\nu^{-\mu}(\cos{\theta}) &=& -\frac{\pi}{2\nu^\mu}\left(\frac{\theta}{\sin{\theta}}\right)^\frac{1}{2}\left[Y_\mu\left(\left(\nu+\frac{1}{2}\right)\theta\right)+{\mathcal O}\left(\frac{1}{\nu}\right){\rm env}Y_\mu\left(\left(\nu+\frac{1}{2}\right)\theta\right)\right], \nonumber
 \ea
 for $0<\theta<\pi-\delta$ with $\delta$ fixed and $\nu\rightarrow\infty$. The envelope functions are treated in \cite[Sec. 2.3.1]{cohl}.  

As may be seen through a comparison with Eqs.\ (\ref{Dcut3}) and (\ref{Ccut3}), the results of the two approaches agree for $\nu\gg 1$, with the simple approximations given here in Theorems 2-4 applying in sectors in $z=\cos{\theta}$ for $0\leq \theta\leq\pi$, and the uniform results holding for for $\theta$ bounded away from $\pi$. The Bessel function expansions derived here also reproduce the first $n$ powers of $(1-z)$ in the Legendre functions properly for $z\rightarrow 1$ when the Bessel functions through order $\mu+n$ are included.

The pre-factors and the variable in the uniform approximations are, unfortunately, awkward for physical applications to scattering theory, where, {\em e.g.},
\be
 \sqrt{2j(j+1)(1-\cos{\theta})}=\sqrt{j(j+1)q^2/p^2}=qb \nonumber
 \ee
rather than $(j+\frac{1}{2})\theta$ is the natural variable. Here $j=\nu$ is conserved angular momentum in the scattering, $q$ is the invariant momentum transfer, $p$ is the momentum of the particles in the center-of-mass system, and $b$ the impact parameter or point of closest approach in the free Schr{\"o}dinger equation.  The pre-factors also disrupt the useful connection between partial-wave series in Legendre functions and Fourier-Bessel transforms in the theory of particle scattering; see, for example, \cite[Appendix B]{DurandChiu}  These problems not encountered with the expansions derived here, Eqs.\ (\ref{Bessel_rep}) and (\ref{D_Bessel2}) for $\theta\approx 0$, and Eqs.\ \ref{DcutY}) and (\ref{CcutY}) for $\theta\approx\pi$. 

Cohl, Dang, and Dunster \cite[Secs. 2.3.1, 2.4.1]{cohl} also treat the limits $\nu\rightarrow\infty$ and $\nu\rightarrow\pm i\infty$ for $z=\cosh{\theta}\in(1,\infty)$ using uniform expansions. Their results in terms of Bessel functions agree in form and error estimate with the simple asymptotic expressions in Theorems 1 and 2  for $\nu\theta\gg 1$, but also extend smoothly to $\theta=0$, $z=1$, the region treated separately in the Bessel function expansions derived here. They do not treat the more complicated cases of complex $z$ and $\nu$, to which the results of Theorems 1 and 2 results apply directly, again away from $z=\pm 1$.\\

Acknowledgments: The author  would  like to thank the Aspen Center for Physics for its hospitality and for its partial support of this work under NSF Grant No. 1066293. He would also like to thank Dr.\ Howard Cohl for raising the questions that led to this work.

\bibliographystyle{unsrt}
\bibliography{Gegenbauer_functions.bib}

\end{document}